%
%

%
%


\magnification 1200
\input amstex
\documentstyle{amsppt}
\NoBlackBoxes
\NoRunningHeads
\TagsOnRight


\define\vre{\varepsilon}
\define\hd{Hausdorff dimension}
\define\hs{homogeneous space}
\define\df{\overset\text{def}\to=}
\define\un#1#2{\underset\text{#1}\to#2}
\define\br{\Bbb R}
\define\bn{\Bbb N}
\define\bz{\Bbb Z}
\define\bq{\Bbb Q}
\define\bc{\Bbb C}
\define\ba{badly approximable}

\define\di{Diophantine}
\define\da{Diophantine approximation}
\define\de{Diophantine exponent}

\define\vx{\bold x}
\define\vy{\bold y}

\define\vv{\bold v}

\define\vq{\bold q}

\define\vc{\bold c}
\define\vf{\bold f}
\define\vt{\bold t}


\define\nz{\smallsetminus \{0\}}

\define\cag{$(C,\alpha)$-good}

\define\vwa{very well approximable}
\define\supp{\operatorname{supp}}
\define\const{\operatorname{const}}

\topmatter
\title 
Diophantine exponents of measures: a dynamical approach
\endtitle  

\author { Dmitry Kleinbock} \\ 
  { \rm 
   Brandeis University} 
\endauthor

    \address{ Dmitry Kleinbock,  Department of
Mathematics, Brandeis University, Waltham, MA 02454-9110}
  \endaddress

\email kleinboc\@brandeis.edu \endemail

  \thanks 
Supported in part by NSF
Grant DMS-0239463.
\endthanks       


\abstract 
We place the theory of  metric \da\ on manifolds
into a broader context of studying \di\ properties of points
generic with respect to certain measures on $\br^n$. 
The correspondence 
 between multidimensional \da\  and dynamics of lattices in 
Euclidean spaces is discussed in an elementary way, and several recent results 
obtained
by means of this correspondence are surveyed.

 \endabstract

\date {March  2004} \enddate


\endtopmatter
\document

\heading{1. Introduction}
\endheading 

We start by recalling several basic facts from the theory of 
 \da.  
For $v>0$ and 
$n\in\bn$,  say that
$\vy = (y_1,\dots,y_n)\in\br^n$ is {\sl $v$-approximable\/} (notation: $\vy\in\Cal
W_v$) if  there are infinitely many 
$\vq= (q_1,\dots,q_n)\in
\bz^n$ such that
$$
 |y_1q_1 + \dots+ y_nq_n + p|   < \|\vq\|^{-v}  \tag 1.1
$$
for some
$p\in\bz$. It will be convenient   to interprete points of $\vy\in\br^n$ as
{\it row vectors\/} (or linear forms) and integers
$\vq\in
\bz^n$ as {\it column vectors\/}, denoting both by boldface letters. This way $y_1q_1 +
\dots + y_nq_n$ can be written as $\vy\vq$, and (1.1) as 
$$
 |\vy\vq + p|   < \|\vq\|^{-v}  \,. \tag 1.2
$$
Hopefully such notation will
cause no confusion.

Then let us define the {\sl \de\/} $\omega(\vy)$ of $\vy$ by
$$
\omega(\vy) \df \sup\{v\mid \vy\text{ is $v$-approximable}\}\,.
$$
Note that the above definition, unlike the previous one, does not depend on the choice 
of the norm $\|\cdot\|$. We will however always work with the sup-norm, $\|\vx\| =
\max_i|x_i|$.

 It is well-known and easy to see that one has $n \le \omega(\vy) \le
\infty$ for all
$\vy$, and $ \omega(\vy) = n$ for $\lambda$-almost every $\vy\in\br^n$,
where $\lambda$ stands for Lebesgue measure on $\br^n$. Vectors $\vy$
with $ \omega(\vy) > n$ are usually called 
\vwa\ (VWA). 

Let us now extend the notion of \de s to measures. Namely, if $\mu$ is a 
locally finite
Borel measure on $\br^n$, let us define the {\sl \de\/} $\omega(\mu)$ of $\mu$ to be the
$\mu$-essential supremum of the function $\vy\mapsto\omega(\vy)$.
In other words, 
$$
\omega(\mu)\df \sup \big\{\,v\bigm| \mu(\{\vy\mid \omega(\vy) > v \}) > 0\big\} = \sup
\big\{\,v\bigm| \mu(\Cal
W_v ) > 0\big\}\,.\tag 1.3
$$
Clearly it only depends on the measure class of $\mu$.

Very often (equivalence classes of) measures $\mu$ that we are going to consider will be
naturally associated with subsets $\Cal M$ of $\br^n$ supporting $\mu$. For example,
if $\Cal M$ is a smooth submanifold of $\br^n$, we will be taking (the class of) $\mu$
to be (that of) the Riemannian volume on $\Cal M$, that is, the pushforward
$\vf_*\lambda$ of
$\lambda$ by any smooth map $\vf$ parametrizing
$\Cal M$. In this case we will define the \de\ $\omega(\Cal M)$ of $\Cal M$ to be equal
to that of
$\mu$. From what was said it follows that
$\omega(\mu)
\ge n$ for any
$\mu$, and $ \omega(\lambda)= \omega(\br^n)$ is equal to $n$. 

This justifies the terminology which has been introduced to \da\ on manifolds
by V.\ Sprind\v zuk: let us say that a measure $\mu$ on $\br^n$ (resp.\ a subset
$\Cal M$ of $\br^n$) is  {\sl
extremal\/} if
$\omega(\mu)$ (resp.\ 
$\omega(\Cal M)$) is equal to $n$, that is, 
attains the smallest
possible value; equivalently, if $\mu$-a.e.\ $\vy\in\br^n$ is not VWA.

\smallskip

The theory started  
with considering the  map
$$
\vf(x) = (x,x^2,\dots,x^n)
\,.\tag 1.4
$$ 
The extremality of $\vf_*\lambda$ for $\vf$ as above
was conjectured in 1932 by K.~Mahler \cite{M}
and proved in 1964 by Sprind\v zuk  \cite{Sp1, Sp2}. In about the same time
W.\ Schmidt \cite{Sc1} proved the extremality of $\vf_*\lambda$
when $\vf:I\to \br^2$, $I\subset \br$, is
$C^3$ and satisfies $$\left|\matrix
f_1'(x) & f_2'(x)\\ f_1''(x) & f_2''(x)
\endmatrix\right| \ne 0\quad \text{for $\lambda$-a.e. }x\in I\,.$$ 
Since then, a lot of attention has been devoted to showing that measures $\vf_*\lambda$
are extremal for other smooth maps $\vf$. 

\smallskip

To describe a broad class of examples, let us recall the following definition. 
Let
$U$ be an open subset of
 $\br^d$   and let  $\vf = (f_1,\dots,f_n)$ be a
  $C^k$ 
map $U\to \br^n$. For $l\le k$ and $\vx\in U$, say that 
$\vf$  is {\sl $l$-nondegenerate  at $\vx$\/} if 
$$
\aligned
\text{$\br^n$ is spanned by partial derivatives of $\vf$ at
$\vx$ of order up to $l$}\,.
\endaligned
\tag 1.5
$$
We will say that
$\vf$ is 
 {\sl nondegenerate  at $\vx$\/}  if (1.5) holds for
some $l$. 
If $\Cal M$ is a $d$-dimensional 
submanifold of $\br^n$, we will say
that $\Cal M$ is {\sl nondegenerate   at $\vy\in \Cal M$} if any
(equivalently, some)  diffeomorphism
$\vf$ between an open subset $U$  of $\br^d$ and a neighborhood of
$\vy$ in $\Cal M$ is   nondegenerate   at $\vf^{-1}(\vy)$. 

In 1996  \cite{KM1}  the following theorem was proved, generalizing the aforementioned
results of Sprind\v zuk and Schmidt:

\proclaim{Theorem 1.1} Let $\Cal M$ be a  smooth 
$d$-dimensional submanifold of $\br^n$ which is nondegenerate at its almost every point.
Then $\Cal M$ is extremal. Or, slightly more generally, if $U$ is an open subset of 
$\br^d$ and $\vf:U\to
\br^n$ is nondegenerate  at $\lambda$-almost every point of $U$, then $\vf_*\lambda$ is
extremal. 
\endproclaim

Note that a real analytic version of Theorem 1.1 was conjectured by Sprind\v zuk
\cite{Sp3} in 1980. The case $n = 3$, $d = 1$ was established earlier by
Beresnevich and Bernik \cite{BB}, and  later an alternative proof of the general case
was found by Beresnevich \cite{Be}.

The method of proof in  \cite{KM1} was dynamical in nature, and the main purpose of this
note is to explain its main ideas and a possibility to use them to
solve more general problems. Speaking of which, it
seems natural not to restrict oneself to smooth measures on submanifolds, and thus
ask

\proclaim{Question 1.2} What other measures  on $\br^n$  can be shown to be extremal?
\endproclaim

In fact, pushing it even further, one can ask

\proclaim{Question 1.3} For what other measures $\,\mu$ on $\br^n$  can one compute
or estimate
$\omega(\mu)$?
\endproclaim

In the present paper we attempt to provide  partial answers to both questions. 
Note that the set-up naturally generalizes to
so-called {\it Khintchine-type theorems\/}, where one replaces the right hand side 
of (1.1) by a function of $\|\vq\|$, or, more generally, of $\vq$. Many  results
for smooth measures on manifolds have been obtained in recent years, but those
will be outside of the scope of this paper.

The goals of the paper are:

\roster
\item"$\bullet$" to describe a correspondence between \da\ and dynamics;
\item"$\bullet$" to state the main ``quantitative nondivergence'' estimate, 
which serves as the main tool for proving many results in both homogeneous dynamics
and \da;
\item"$\bullet$" to list results one can derive using this method, some of them old
and some new and yet unpublished.
\endroster

This was roughly the outline of the author's talk given at the conference on
``Diophantine analysis, uniform distributions and applications'' 
in  Minsk, Belarus in August 2003. The hospitality of the organizers of this
conference is gratefully acknowledged. Thanks are also due to Victor Beresnevich
and Barak Weiss for useful comments.

\heading{2. \da\ and dynamics}
\endheading 

The correspondence between \di\ properties of vectors in $\br^n$ and dynamical
properties of lattices in $\br^{n+1}$ dates back for $n = 1$ to E.\ Artin,
and for $n > 1$ to the work of Schmidt \cite{Sc2} and Dani \cite{D}.
Here we present a condensed exposition of the main principle behind a reduction
of Theorem 1.1 to a dynamical statement. Note that a similar exposition can
be found in survey papers \cite{K1} and \cite{Ma2}, as well as in Chapter IV
of \cite{St}.

We are going to pick $v_0 \ge n$ and $\vy\in\br^n$ with $\omega(\vy)  > v_0$.
According to the definition, this means that for some $v > v_0$ inequality (1.2) is
satisfied for infinitely many $(p,\vq)$. Equivalently, for infinitely
many $s\in\bn$ the following system has an integer solution $(p,\vq)$:
$$
\cases|\vy\vq + p|   < 2^{-vs} \\2^s \le \|\vq\| < 2^{s+1}\,. \endcases
$$
We will drop the first of the inequalities in the second line,
and  conclude that $\omega(\vy)  > v_0$ implies that for some $v > v_0$ there exist
infinitely many $s\in\bn$ such that the  system
$$
\cases|\vy\vq + p|   < 2^{-vs} \\ \|\vq\| < 2^{s+1}\endcases\tag 2.1
$$
has an integer solution $(p,\vq)$ with $\vq\ne 0$.

The latter system can be conveniently written in a matrix form. Namely, $\vy$ gives rise
to
 $$\ 
u_\vy \df \left(\matrix
1 & \vy  \\
0 & I_n
\endmatrix \right)\,,
$$  
and 
the right hand sides of the inequalities in (2.1) define a certain rectangular box 
in $\br^{n+1}$, namely
$$
B_{v,s} = \big\{(x_0,x_1,\dots,x_n)\bigm| |x_0| < 2^{-vs},\ |x_i| < 2^{s+1}\text{ for }i
= 1,\dots,n\big\}\,.
$$
Thus
$\omega(\vy)
> v_0$ implies that for some $v > v_0$ there exist infinitely many $s\in\bn$ 
such that 
$u_\vy\bz^{n+1}\cap 
B_{v,s}\ne\{0\}$. 

\smallskip

The next step is to transform this box into a cube, which is naturally easier
to work with than a thin and flat rectangular box $
B_{v,s}$. This is where dynamics comes into
play. There is only one way to undertake such a transformation preserving the volume:
one needs to act by 
$$
g_t = \text{\rm diag}(2^{nt},
2^{-t},\dots,2^{-t})
\,.\tag 2.2
$$
Note that the volume of $
B_{v,s}$ is equal to $2^{2n + 1 - (v - n)s}$, which tends to $0$ when $s\to\infty$
(recall that $v$ is chosen to be bigger than $v_0 \ge n$). If for any $s$ one chooses
$t > 0$ such that   $
g_tB_{v,s}$ is a cube, then it is clear that the sidelength of this cube will be very 
small for large $s$. In fact, an elementary computation shows that $s$ and $t$ 
in this situation are not far from each other, and 
the sidelength of $
g_tB_{v,s}$ is equal to $C_n2^{-\gamma t}$, where $C_n$ is an explicit constant
depending only on $n$ and 
$$
\gamma = \gamma(v) = \frac{v-n}{n(v+1)}\,.
$$
We have almost proved

\proclaim{Proposition  2.1} Suppose that $\omega(\vy)> v_0$ for some $\vy\in\br^n$
and  $v_0 \ge  n$. Then for some $\gamma > \gamma(v_0)$ there exist infinitely many
$\,t\in\bn$ such that 
$$
\text{the lattice $g_{t}u_\vy\bz^{n+1}$ has a nonzero vector
of norm less than }2^{-\gamma t}\,.\tag 2.3$$
\endproclaim

\demo{Proof} It remains to observe that taking $\gamma$ between 
$\gamma(v)$ and  $\gamma(v_0)$ allows one to get rid of the constant $C_n$, as well as
to replace  every $t$ chosen as above by its integer part.
\qed\enddemo

In fact, the converse to this proposition is also true and easy to prove, see 
\cite{K2} or \cite{KM2}, 
but it will not be needed here.

\smallskip

It might be  helpful
for the understanding to discuss the
geometric meaning of the conclusion of the above proposition. Denote by 
$\Omega_k$ the space of lattices in $\br^k$ of
covolume $1$, and let
$$
\Omega_k(\vre) \df \big\{\Lambda \in\Omega_k \bigm| \|\vv\| < \vre\text{ for some
}\vv\in\Lambda\nz\big\}\,.
$$
Then (2.3) can be  written as $g_{t}u_\vy\bz^{n+1}\in \Omega_{n+1}(2^{-\gamma
t})$.
It is well known that $\Omega_k$ is noncompact, but the
complement of $\Omega_k(\vre)$ in $\Omega_k$ is compact for
any positive $\vre$, and, further, any bounded subset of $\Omega_k$ belongs to
such a complement for some $\vre > 0$ (Mahler's Compactness
Criterion, see \cite{R} or \cite{BM}). 
Thus vectors with large \de s
give rise to
$g_t$-trajectories in the space of lattices with ``fast enough growth''. See
\cite{K1, K2} for more details.

Here is an application of Proposition  2.1:

\proclaim{Corollary  2.2}  Let $U$ be an open subset of
$\br^d$, $\mu$ a measure on $U$,   and let $\vf$ be a  map from $U$ to $\br^n$.
Take $v \ge n$, and suppose that for $\mu$-a.e.\ $\vx_0\in
U$ one can find a ball $B\subset U$ centered in $\vx_0$ such that for any $\gamma >
\gamma(v)$ one has
$$\sum_{t = 1}^\infty\mu\big(\big\{\vx\in B \bigm| g_{t}u_{\vf(\vx)}\bz^{n+1}\in
\Omega_{n+1}(2^{-\gamma t})\big\}\big) < \infty\,.\tag 2.4
$$
Then $\omega(f_*\mu) \le v$.
\endproclaim

\demo{Proof} Indeed, in view of the
Borel-Cantelli Lemma,  (2.4) implies that 
$$\mu\big(\big\{\vx\in B \bigm| g_{t}u_{\vf(\vx)}\bz^{n+1}\in \Omega_{n+1}(2^{-\gamma
t})\text{ for infinitely many }t\big\}\big) = 0\,.
$$
Hence it follows from the assumption and Proposition  2.1 that for $\mu$-a.e.\ $\vx_0\in
U$ one can find a ball $B$ centered in $\vx_0$ such that
$\mu\big(\vx\in B \mid \omega\big(\vf(\vx)\big)> v\}\big) = 0
$, and the latter, in view of the definition (1.3), implies that $\omega(\mu) \le v$.
\qed\enddemo

Summarizing the above discussion, we can observe that an upper estimate for
$\omega(f_*\mu)
$, and in particular the extremality of $f_*\mu$, can be derived from knowing that
sets of the form
$$
\big\{\vx\in B \bigm| g_{t}u_{\vf(\vx)}\bz^{n+1}\in
\Omega_{n+1}(\vre)\big\}\tag 2.5
$$
have small enough measure. In other words, the $g_t$-translate of the pushforward
of $\mu$ by the map $\vx\mapsto u_{\vf(\vx)}\bz^{n+1}$, $B\to \Omega_{n+1}$, does
not assign too much weight to the ``neighborhood of infinity''  $\Omega_{n+1}(\vre)$ in
the space of lattices. The latter is precisely a consequence of so-called ``quantitative
nondivergence estimates'', to be discussed in the next section.

\heading{3. Quantitative nondivergence and proof of Theorem 1.1}
\endheading 

In order to state a general result which can be used to estimate the measure of 
sets (2.5), we need to introduce some notation and definitions.
For a ball $B=B(\vx, r)$ in $\br^n$ and $a > 0$, we denote $B(\vx,ar)$ by $a
     B$. If  $B$ is a ball in $\br^n$ and $f$ is a real-valued function on
$\br^n$, let
$$\Vert f \Vert_{B} \df\sup_{\vx \in B} |f(\vx)|\,;
$$
and if $\mu$  is  a measure on $\br^n$ such that $\mu(B) > 0$, we define
$\Vert f
\Vert_{\mu,B}$  to be equal to $\Vert f \Vert_{B\,\cap\, \supp\, \mu}\,$.

Given $D \ge 1$, say that a  measure
$\mu$ on $\br^n$ 
 is 
{\sl $D$-Federer on\/}  $U$ if
$${\mu\left(3  B \right)}\leq D{\mu\left(B\right)}
$$
for every ball $B$ centered in $\supp \, \mu$ with $3B \subset U$.
We will say that a measure is  {\sl 
Federer\/} if  for $\mu$-a.e.\ point of
$\br^n$  there exist a neighborhood $U$ of this point and  $D > 0$ such that $\mu$ is
 $D$-Federer on $U$.

Clearly $\lambda$ and, more generally, volume measures on smooth submanifolds 
satisfy the above condition. But many other natural measures can also be proved
to be  Federer. See  \cite{KLW, KW, MU, S} for
examples.

\smallskip

The next definition involves a very important property of certain functions
$f$ with respect to certain measures $\mu$.
Given $C$, $\alpha > 0$,  a subset $U$ of $\br^n$,
a
measure
$\mu$ on
$U$
        and a
real-valued function $f$ on $U$, say that $f$ is
        {\sl $(C,\alpha)$-good
on  $U$
with respect to
$\mu$\/}
        if for 
any open ball
$B 
\subset U
$ centered in $\supp\,\mu$
and any
$\vre > 0$ one has
$${\mu\big(\big\{x\in B \bigm| |f(x)| < \vre\big\}\big)} \le C
\left(\frac{\varepsilon}{\Vert f\Vert_{\mu, B}}\right)^\alpha{\mu(B)}\,.
$$
The primary example is given by polynomial maps. See  \cite{KM1, BKM, KLW}, 
as well as \S
4 of the present paper, for various other examples.

\smallskip

We are now ready to state our main estimate. It was proved in 
 \cite{KM1} (Theorem 5.4) in the case $\mu = \lambda$ , and then generalized 
in \cite{KLW} and \cite{KT}.

 \proclaim{Theorem 3.1}  For  $d,n\in \bn$, let a ball  $B\subset \br^d$,  a
measure
$\mu$ on
$\br^d$ such that $B$ is   centered at $\supp\,\mu$ and $\mu$ is $D$-Federer
on $\tilde B \df 3^{n+1}{B}$,  and a continuous map  $\vf = (f_1,\dots,f_n):\tilde B \to
\br^n$ be given. Suppose also that for some $C,\alpha > 0$ and $0 < \varrho <
\frac1{n+1}$ the following two conditions hold:
\roster
\item"(i)" 
for any $\vc = (c_0,c_1,\dots,c_n) \in \br^{n+1}$, the function $c_0 + \sum_{i =
1}^nc_if_i$ is
\cag\ on $\tilde B$ w.r.t.\ $\mu$; 
\item"(ii)" 
for any $\vc\in \br^{n+1}$ with $\|\vc\| \ge 1$, $\|c_0 + \sum_{i =
1}^nc_if_i\|_{\mu,B} \ge \varrho$.
\endroster
Then 
 for any  positive $ \vre \le \varrho$ and any $t > 0$ one has
$$
\mu\big(\big\{\vy\in B \bigm| g_{t}u_\vy\bz^{n+1}\in
\Omega_{n+1}(\vre)\big\}\big) \le  ({n+1})C\big(N_dD^2\big)^{n+1} \left(\frac\vre \varrho
\right)^\alpha  \mu(B)\,,
$$
where $N_d$ (the so-called Besicovitch constant 
of $\br^d$) depends  only on $d$.
\endproclaim

The proof of this theorem is not easy. But fortunately most of it has been around
since the early 1970s, when Margulis \cite{Ma1} proved that unipotent flows on the space
of lattices do not go to infinity. In fact, his proof applies verbatim if one replaces
a unipotent subgroup by a polynomial map. Later it was realized that the only
way the polynomiality of the map is used is via the \cag\
property, and that it can also produce a quantitative strengthening of non-divergence
to infinity, namely an estimate for a measure of the intersection of the trajectory
and  a small ``neighborhood of infinity''  in the space of lattices.

\smallskip

Let us now show how this theorem can be applied to \da. Let  an open subset $U$ of
$\br^d$, 
a measure $\mu$ 
on $U$, $\vx_0\in U$, and  a map  $\vf$ from $U$ to $
\br^n$ be given.

For brevity, let us say that $\vf= (f_1,\dots,f_n)$ is 

\roster
\item"$\bullet$" {\sl
$\mu$-good  at 
$\vx$\/}  if 
$$
\aligned
\text{there exists a neighborhood
$V\subset U$ of $\vx$ and positive $C,\alpha$
such that}\\
     \text{any linear combination of
$1,f_1,\dots,f_n$  is
$(C,\alpha)$-good on
$V$ w.r.t.} \ \mu;\endaligned\tag 3.1
$$

\item"$\bullet$" {\sl $\mu$-nonplanar at
$\vx$\/}
    if
$$
\aligned
\text{for any  neighborhood
$B$ of $\vx$}\\\text{ the
restrictions of $1,f_1,\dots,f_{n}$ to
}B\,\cap &\,\supp\,\mu\\
     \text{
are linearly independent over }&\br;\endaligned\tag 3.2
$$ 
in other words, if $\vf(B \,\cap\,
\supp\,\mu)$  is not contained in any proper affine subspace  of $\br^n$.
\endroster

We can now prove

 \proclaim{Theorem 3.2}  Let $\mu$ be a Federer measure on
$\br^d$, $U$ an open subset of
$\br^d$, and $\vf:U\to\br^n$ a continuous map which is $\mu$-good and
$\mu$-nonplanar  at $\mu$-almost every point of $U$.  Then $\vf_*\mu$
is extremal.
\endproclaim

\demo{Proof} Take  $\vx\in U\,\cap \,\supp\,\mu$ satisfying (3.1) and (3.2), 
and let $B$ be a ball centered at $\vx$ such that $3^{n+1}{B}$ is contained in $V$.
Then condition (i) of Theorem 3.1 will be satisfied for some positive $C,\alpha$, and 
the existence of $\varrho > 0$ satisfying (ii) follows from the compactness of the
unit sphere in
$\br^{n+1}$. Thus Theorem 3.1 applies and it follows for any small enough positive $\vre$
 and any $t > 0$ one has
$$\mu\big(\big\{\vx\in B \bigm| g_{t}u_{\vf(\vx)}\bz^{n+1}\in
\Omega_{n+1}(\vre)\big\}\big) \le \const\cdot \vre ^\alpha \,,
$$
with the constant independent of $\vre$ or $t$. Putting $\vre = e^{-\gamma t}$
for an arbitrarily small $\gamma > 0$ and using Corollary 2.2 finishes the proof.
\qed\enddemo

The above theorem first appeared in \cite{KLW}  in a slightly disguised version:
there $n$ was equal to $d$, $\vf$ was the identity map, and  
 conditions sufficient for extremality were stated in terms of $\mu$. But the proof
given there (which itself is a generalization of the argument from \cite{KM1})
in fact easily yields the result stated above. The argument was generalized
even further in
\cite{KT}, where the wording was similar to that of the present paper.

Note also that by definition, the set
$$
\{\vx\mid \vf \text{ is $\mu$-nonplanar at }\vx\}
$$
is closed, so the nonplanarity  of $\vf$ at $\mu$-almost every point is equivalent to
the same at  every $\vx\in U\,\cap \,\supp\,\mu$.

In order to see that Theorem 1.1 is a special case of Theorem 3.2, it suffices 
to show that a smooth map  $\vf:U\to\br^n$ 
is $\lambda$-good and
$\lambda$-nonplanar  at every point where it is nondegenerate. The nonplanarity is
straightforward (indeed, the nondegeneracy of $\vf$  at
$\vx$ clearly implies the existence of a neighborhood $B\ni \vx$ such that  $\vf(B
)$  is not contained in any proper affine subspace  of $\br^n$). And the \cag\ property
of linear combinations of $1,f_1,\dots,f_n$ basically follows from the fact that
locally $\vf$ can be approximated by a polynomial map, and is proved in \cite{KM1}.

\smallskip
Now it seems to be worthwhile to compare the method of proof of Theorem 1.1 discussed
above with the standard approach based on Sprind\v zuk's solution of Mahler's problem
and carried out in \cite{Be}.
As we have seen, the correspondence between \da\ and dynamics is quite natural and easy
to explain. Also, the crucial measure estimate (Theorem 3.1) is the only hard part of
the argument, the rest is relatively easy. Another advantage is a chance to
work with non-smooth objects   --  we will mention in the   next section how Theorem 3.2
gives rise to a wide variety of examples of  extremal measures which are not volume
measures on smooth submanifolds.
 Further,
as will also be discussed in the   next section, the dynamical approach can be perturbed
in many directions and allows many generalizations and modifications of the 
result proved above. 

However the standard methods have a number of obvious advantages as well. The dynamical
approach is hard to use when more precise results are needed, for example when the
goal is to prove the divergence case of Khintchine-type theorems,
or compute/estimate the \hd\ of the set of $v$-approximable points on a manifold. 
See \cite{Be, BDV1, BDV2} for examples of such results.
Roughly speaking, the correspondence between approximation and dynamics is powerful but
coarse, so that a substantial amount of information is being lost in transmission.

\heading{4. Beyond Theorem 1.1}
\endheading

\subheading{4.1} As was mentioned before, one of the main advantages of the method is
that the assumptions of Theorem 3.2 are much less restrictive that those of Theorem 1.1.
Here is an example. Following \cite{KLW}, given $C,\alpha > 0$ and $U\subset \br^d$, say
that $\mu$ is
{\sl absolutely $(C, \alpha)$-decaying on\/}  $U$
if
for any non-empty open ball $B \subset U$
of radius
$r$
centered
in $\supp\,\mu$,
any affine hyperplane $\Cal L \subset
\br^n$, and any $\varepsilon >0$
one has
$${\mu \big( B \cap \Cal L^{(\varepsilon)} \big) }
\le C
\left(
\frac{\varepsilon}{r
}
\right)^{\alpha}{\mu(B)}\,,
$$
where $\Cal L^{(\varepsilon)}$ stands for the $\vre$-neighborhood of $\Cal L$.
We will say that a measure is {\sl absolutely decaying\/}  if  for $\mu$-a.e.\ point of
$\br^d$  there exist a neighborhood $U$ of this point and $C,\alpha >
0$  such that $\mu$ is absolutely
$(C,\alpha)$-decaying on $U$. Measures which are absolutely decaying and Federer 
were called {\sl absolutely friendly\/} in \cite{PV}, see \cite{KLW} for justification
of this terminology. 

If $\mu$ is absolutely
decaying, it easily follows  that $$\mu (\Cal
L
) = 0\text{ for any affine hyperplane }\Cal L\subset
\br^n\tag 4.1
$$
 (i.e., in the terminology introduced
in \S 3, the identity map $\br^n\to \br^n$ is $\mu$-nonplanar at every point of
$\supp\,\mu$). It also follows that the identity map is $\mu$-good at $\mu$-almost every
point.  In \cite{KLW}, much more than that has been proved:

 \proclaim{Proposition 4.1}  Let $\mu$ be an absolutely friendly measure on
$U\subset\br^d$, and let $\vf:U\to\br^n$ be a $C^{l+1}$  map  which is
$l$-nondegenerate at  $\vx_0\in U$.  Then $\vf$
is  {\rm (a)} $\mu$-good and {\rm (b)}
$\mu$-nonplanar  at $\vx_0$.
\endproclaim

Part (b) is straightforward from (4.1), but part (a) is nontrivial and can be thought 
of as
a generalization of the case $\mu = \lambda$  worked out in \cite{KM1}.
However note that the method of proof is completely different, has an advantage of
producing a better esponent $\alpha$ in many cases, but has a slight
disadvantage of requiring an extra derivative.

The above proposition immediately implies that  $\vf_*\mu$
is extremal whenever $\mu$ is absolutely  friendly and $\vf$ is 
nondegenerate. The result is interesting  even when $\vf$ is the identity map
(recently an alternative proof of the latter special case has been worked out in
\cite{PV}). Indeed, one can exhibit a wide variety of absolutely  friendly measures
supported on very peculiar sets, such as self-similar or self-conformal fractals.
The prime example is the middle-third Cantor set ${\Cal C}$ on the real line: its
extremality (or, more precisely, the extremality of the natural measure $\mu_{\Cal C}$
it supports) was established by Barak Weiss in \cite{W}, and later in \cite{KLW} higher
dimensional generalizations of $\mu_{\Cal C}$ were shown to satisfy conditions
sufficient for extremality. More examples have been recently found by M.\ Urba\'nski
\cite{U1, U2}.

\subheading{4.2} The dynamical approach is very useful in handling the so-called
{\it multiplicative\/} generalization of the problems discussed in the introduction.
Namely, define $\Pi_{\sssize +}(\vq) 
\df \prod_{q_i \ne 0 } |q_i|$,  say that
$\vy = (y_1,\dots,y_n)\in\br^n$ is {\sl $v$-multiplicatively approximable\/} (notation:
$\vy\in\Cal W^\times_v$) if  there are infinitely many 
$\vq\in
\bz^n$ such that
$$
 |\vy\vq + p|   < \Pi_{\sssize +}(\vq)^{-v/n}  
$$
for some $p\in\bz$, and then define  {\sl multiplicative \de s\/}:

\roster
\item"$\bullet$" $\omega^\times(\vy)$ of $\vy$ by
$
\omega^\times(\vy) \df \sup\{v\mid \vy\in \Cal
W^\times_v\}\,,
$
\item"$\bullet$" $\omega^\times(\mu)$ of $\mu$ by
$
\omega^\times(\mu)\df 
\sup
\big\{\,v\bigm| \mu(\Cal
W^\times_v ) > 0\big\}\,.
$
\endroster
 
It is  easy to see that  $\omega^\times(\vy)$ is not less than $ \omega(\vy)$ for
all
$\vy$, and yet $ \omega^\times(\vy) = n$ for $\lambda$-a.e.\ $\vy\in\br^n$,
that is, $\omega^\times(\lambda) = n$.  Following  Sprind\v zuk, say that 
$\mu$  is  {\sl  strongly
extremal\/} if
$\omega^\times(\mu) = n$. 
The multiplicative analogue of Theorem 1.1 (more precisely, 
of its real analytic version) was conjectured by Sprind\v
zuk in 1980 (earlier A.\ Baker \cite{B} conjectured that the curve parametrized by $\vf$
as in (1.4) is strongly
extremal) and proved in \cite{KM1}. Likewise, the following can be proved:

 \proclaim{Theorem 4.2}  Let   $\mu$  and $\vf$ be as in 
Theorem 3.2.  Then $\vf_*\mu$
is strongly extremal.
\endproclaim

Thus in all the examples mentioned in \S 4.1,
$\vf_*\mu$  actually happens to be  strongly extremal. The proof of the stronger
statement is based on using the multi-parameter action of 
$$
g_\vt = \text{\rm diag}(2^{t},
2^{-t_1},\dots,2^{-t_n})
\,,\quad\text{where}\quad t = t_1 + \dots + t_n\,,
$$
instead of (2.2).

\subheading{4.3} Obvious examples of non-extremal manifolds are provided
by proper affine subspaces of $\br^n$ whose coefficients are well enough 
approximable by  rational numbers. On the other hand, it is clear from
a Fubini argument that almost all translates of any given subspace are extremal.
In \cite{K2} the method of \cite{KM1} was pushed further to produce criteria for the
extremality, as well as the strong extremality, of an affine subspace $\Cal L$ of
$\br^n$. Further, it was shown that if $\Cal L$ is extremal 
(resp.\ strongly extremal), then so is any   smooth submanifold of $\Cal L$  which is
nondegenerate in
$\Cal L$ at its a.e.\ point. (The latter property is a straightforward generalization
of the definition  of nondegeneracy in $\br^n$: a map $\vf$
  is {\sl nondegenerate in $\Cal L$ at $\vx$\/} if 
the linear part of $\Cal L$ is spanned by partial derivatives of $\vf$ at $\vx$.)
In other words, extremality and strong extremality pass from affine subspaces
to their nondegenerate submanifolds.

A more precise analysis makes it possible to study \de s of measures with supports
contained in proper affine subspaces of
$\br^n$. Namely, in \cite{K4} it is shown how to compute $\omega(\Cal L)$
and $\omega^\times(\Cal L)$ for any $\Cal L$. Moreover, the following generalization
of Theorem 3.2 and its stronger multiplicative form is obtained:

 \proclaim{Theorem 4.3}  Let  $\mu$ be a Federer measure on
$\br^d$, $U$ an open subset of
$\br^d$, $\Cal L$ an affine subspace of
$\br^n$, and let $\vf:U\to\Cal L$ be a continuous map which is $\mu$-good
and
$\mu$-nonplanar in $\Cal L$ at $\mu$-almost every point of $U$.  Then $\omega(\vf_*\mu)
= \omega(\Cal L)$ and $\omega^\times(\vf_*\mu)
= \omega^\times(\Cal L)$.
\endproclaim

Here we say that $\vf$ is {\sl $\mu$-nonplanar in $\Cal L$ at $\vx$\/} if 
for any  neighborhood
$B$ of $\vx$, the $\vf$-image of $B \,\cap\,
\supp\,\mu$ is not contained in any proper affine subspace  of
$\Cal L$,
thus generalizing the definition from \S 3. It is easy to see  that smooth maps
$\vf:U\to\Cal L$ are $\mu$-good
and
$\mu$-nonplanar in $\Cal L$ at every point at which they are nondegenerate in $\Cal L$.

\subheading{4.4} Another application  concerns badly
approximable vectors. Recall that $\vy \in\br^n$ is called {\sl badly approximable\/}
if  there exists $c > 0$ such that for any 
$\vq\in
\bz^n\nz$ and $p\in\bz$ one has 
$$
 |\vy\vq + p|   > c \|\vq\|^{-n}\,.  
$$
Denote the set of \ba\ vectors in $\br^n$ by $\bold{BA}$. It is a theorem of Dani
\cite{D} that
$\vy\in\bold{BA}$  iff the trajectory
$$\{g_{t}u_\vy\bz^{n+1}\mid t > 0\}\,,$$
where $g_{t}$ and $u_\vy$ are as in \S 2, does not intersect $
\Omega_{n+1}(\vre)$ for some $\vre > 0$, i.e.\ is bounded in 
$\Omega_{n+1}$.

Using this and Theorem 3.1 it turns out to be possible to find  badly approximable
vectors inside supports of certain measures on $\br^n$. Here is one way to state
the result of the paper \cite{KW}. Denote by $\dim(K)$ the \hd\ of a subset $K$ of
$\br^n$, and for $\beta > 0$ let us say that a measure $\mu$  on $\br^n$ is {\sl
$\beta$-scaling\/} if there
is a positive
$c>0$ such that for every ball $B$ of radius $r$ one has
$\mu\left(B
\right) \leq c \, r^{\beta}
$.
It is well known (mass distribution principle $+$ Frostman's Lemma) that 
$$
\dim(K) = \sup  \{ \beta \mid K \text{ supports a $\beta$-scaling measure} \}\,.
$$

Now let us define
$$
\dim_{af}(K) \df \sup  \{ \beta \mid K \text{ supports a $\beta$-scaling absolutely
friendly measure} \}\,.
$$
Naturally, it is always not bigger than $\dim(K)$ (but can be much less). The following
is essentially proved in
\cite{KW} (see also
\cite{KTV} where a similar result has been announced):

 \proclaim{Theorem 4.4}  For any   compact subset $K$ of
$\br^n$, one has 
$$
\dim_{af}(K\cap\bold{BA} )\ge \dim_{af}(K)\,.$$
\endproclaim

In particular, if a set $K$ of \hd\ $\beta$ supports an absolutely friendly
$\beta$-scaling measure (and many examples of such sets have been found 
in \cite{KLW, KW, U1}), then 
$\dim(K\cap\bold{BA} ) =
\beta$. Note that this also proves that some sets, like the set of \vwa\ vectors which
has \hd\ $n$, do not support any absolutely friendly
 measures.

\subheading{4.5} In all the problems mentioned above, the ground field $\br$ can be
replaced by $\bq_p$, and in fact several fields can be taken simultaneously, thus
giving rise to the  $S$-arithmetic setting where $S = \{p_1,\dots,p_s\}$ is a finite set
of normalized valuations on $\bq$, including or not including the infinite valuation
(cf.\ \cite{Sp2, Z}).
 The space of lattices in $\br^{n+1}$ is replaced there by the space of
lattices in $\bq_S^{n+1}$, where $\bq_S$ is the product of the fields $\br$ and
$\bq_{p_1},\dots,\bq_{p_s}$. This is the subject of the paper \cite{KT} and its sequel,
currently in preparation. Note that  one can also replace $\bq$ by its finite extension
$K$, and
$\bz$ by the integer points of $K$. See also \cite{K3} where the problem of \da\ in
$\bc^n$ is considered, generalizing Sprind\v zuk's solution \cite{Sp2} of the complex
case of Mahler's Conjecture (this involves studying small values of linear forms with 
coefficients in $\bc$ at real integer points).

\subheading{4.6} Finally, let us mention that a generalization of Theorem 3.1 was used in
\cite{BKM} to  estimate a measure of the set of points $\vx$ for which the system
$$
\cases|\vf(\vx)\vq + p|   < \vre \\|\vf'(\vx)\vq |   < \delta \\ |q_i| <
Q_i,\ i = 1,\dots,n
\endcases
$$
has a nonzero integer solution. Here $u_{\vf(\vx)}$ has to be replaced by the matrix
 $$\ 
\left(\matrix
1 & 0 &\vf(\vx)  \\
0 & 1 & \vf'(\vx)  \\0 & 0 & I_n
\endmatrix \right)\,,
$$  
and therefore conditions (i) and (ii) of Theorem 3.1 are replaced by more complicated
conditions, which nevertheless can be checked when $\vf$ is a smooth nondegenerate
map.  This resulted in proving the convergence case of Khintchine-Groshev Theorem 
for nondegenerate manifolds, in both standard  and multiplicative versions. The estimate
was also used in \cite{BBKM} for the proof of the divergence case.


\bigskip
\Refs

\widestnumber\key{BBKM}


\ref\key {B}\by A. Baker  \book Transcendental number theory 
\publ Cambridge Univ. Press \publaddr Cambridge
\yr 1975 \endref

\ref\key Be \by V. Beresnevich  \paper A Groshev type theorem for convergence on
manifolds  \jour Acta Math. Hungar. \vol 94 \yr 2002\pages 99--130 \endref




\ref\key {BB}\by V. Beresnevich and V. Bernik \paper On a metrical theorem of
W.~Schmidt \jour Acta Arith. \vol 75 \pages 219--233 \yr 1996\endref%


\ref\key BBKM \by V. Beresnevich, V. Bernik, D. Kleinbock, and G.\,A.
Margulis 
 \paper Metric Diophantine approximation: the Khintchine--Groshev theorem
for non-degenerate manifolds \jour Moscow Math. J. \yr
2002 \vol 2 \issue 2 \pages 203--225\endref


\ref\key BD \by V. Bernik and
M.\,M. Dodson \book Metric \da\
on manifolds \publ Cambridge Univ. Press \publaddr Cambridge
\yr 1999 \endref

\ref\key BDV1 \by   V. Beresnevich, H. Dickinson and S. Velani
 \paper  Measure theoretic laws for lim sup sets\paperinfo Preprint \yr  2003 \endref Ê

\ref\key BDV2 \bysame \paper  Diophantine approximation on planar curves and the distribution of rational points
\paperinfo Preprint \yr  2003 \endref Ê


\ref\key BKM \by V. Bernik, D. Kleinbock, and G.\,A. Margulis \paper
Khintchine-type theorems  on
manifolds:  the convergence case for standard  and multiplicative
  versions \jour Internat. Math. Res. Notices \yr 2001   
\pages 453--486 \issue 9
\endref


\ref\key BM  \by B. Bekka and M. Mayer 
\book Ergodic theory and topological dynamics of group actions 
on
  homogeneous spaces  
\publ Cambridge University Press \publaddr Cambridge \yr 2000 \endref



\ref\key {D}\by S.\,G. Dani \paper Divergent trajectories of flows on
\hs s and Diophantine approximation\jour
J. Reine Angew. Math.\vol 359\pages 55--89\yr 1985\endref







\ref\key {K1} \by D. Kleinbock \paper Some applications of
homogeneous dynamics to number theory \inbook in: Smooth Ergodic
Theory and Its Applications (Seattle, WA, 1999) \pages 639--660
\bookinfo Proc. Symp. Pure Math. \vol 68  \publ Amer. Math. Soc.
\publaddr Providence, RI \yr 2001 \endref

\ref\key K2 \bysame \paper Extremal subspaces and their
submanifolds \jour Geom. Funct. Anal. \vol 13 \yr 2003 \issue 2 \pages
437--466 \endref

\ref\key K3 \bysame \paper Baker-Sprind\v zuk conjectures for complex analytic manifolds
\inbook in: Algebraic groups and Arithmetic, TIFR, India \toappear \endref


\ref\key {K4} \bysame  \paper  An extension of quantitative nondivergence 
and applications to 
Diophantine exponents \paperinfo in preparation \endref


\ref\key {KLW} \by D. Kleinbock, E. Lindenstrauss, and B. Weiss \paper On fractal
measures and Diophantine approximation \paperinfo Preprint \yr 2003\endref

\ref\key KM1 \by D. Kleinbock and G.\,A. Margulis \paper Flows  on
homogeneous spaces and \da\ on manifolds\jour Ann. Math. \vol 148 \yr
1998 \pages 339--360 
 \endref

\ref\key KM2 \bysame \paper Logarithm laws for flows  on
homogeneous spaces \jour Invent. Math.\vol 138 \pages 451--494 \yr 1999 
\endref




\ref\key KT \by D. Kleinbock and G. Tomanov \paper
Flows on $S$-arithmetic homogeneous spaces and
applications  to metric
Diophantine approximation 
\paperinfo Max Planck Institute Preprint
2003--65 \yr 2003\endref



\ref\key{KTV} \by S. Kristensen, R. Thorn, and S. Velani \paper
Diophantine approximation and \ba\ sets\paperinfo 
in preparation\endref 



\ref\key KW \by D. Kleinbock  and B. Weiss
\paper Badly approximable vectors
on fractals \paperinfo  Preprint\yr  2003 \endref 

\ref\key {M}\by K. Mahler 
\paper \" Uber das Mass der Menge aller $S$-Zahlen \jour Math. Ann. \vol 106 \pages 131--139 \yr 1932\endref

\ref\key {Ma1}\by G.\,A. Margulis
\paper On the action of unipotent group in the space of lattices 
\inbook Proceedings of the Summer School on group representations, (Budapest 1971)\pages
365--370\publ Acad\'emiai Kiado
\publaddr Budapest \yr 1975\endref

\ref\key {Ma2}\bysame \paper Diophantine approximation, lattices and flows on homogeneous
spaces \inbook in: A panorama of number theory or the view from Baker's garden (Z\"
urich, 1999) \pages 280--310 \publ Cambridge Univ. Press \publaddr Cambridge\yr  2002 
\endref



\ref\key {MU} \by D. Mauldin and M. Urba\'nski \paper
Dimensions and measures in infinite iterated function systems
\jour Proc. London Math. Soc. \issue 1 \vol 73 \yr 1996 \pages 105--154\endref


\ref\key PV \by A. Pollington and S. Velani \paper Metric Diophantine approximation and
`absolutely friendly' measures  ÊÊÊÊ\paperinfo Preprint \yr  2003 \endref Ê
\ref \key R \by M.\,S. Raghunathan \book Discrete subgroups of Lie
groups 
\publ Springer-Verlag \publaddr Berlin and New York \yr 1972 \endref%
 

\ref\key S \by H. Sato \paper Global density theorem for a Federer
measure \jour Tohoku Math. J.  \vol 44  \yr 1992 \issue 4 \pages
581--595
\endref


\ref\key {Sc1}\by W. Schmidt \paper Metrische S\"atze \"uber simultane 
Approximation abh\"anginger Gr\"ossen \linebreak \jour Monatsch. Math. \vol 68 \pages
154--166\yr 1964\endref


\ref\key {Sc2}\bysame \paper Diophantine approximation and certain
sequences of 
lattices \jour Acta Arith. \vol 18 \yr 1971 \pages 195--178\endref

 


\ref\key {Sp1}\by V. Sprind\v zuk \paper More on Mahler's conjecture  \jour  Doklady  Akad.  Nauk  SSSR \vol 155 \yr 1964 \pages 54--56  \lang Russian  \transl\nofrills English transl. in  \jour
Soviet Math. Dokl \vol 5 \pages
361--363\yr 1964\endref



\ref\key {Sp2}\bysame \book Mahler's problem in metric number theory \bookinfo Translations of Mathematical
Monographs, vol. 25 \publ Amer. Math. Soc.\publaddr Providence, RI \yr 1969 \endref


\ref\key {Sp3}\bysame  \paper Achievements and problems in
Diophantine approximation theory \jour Russian Math. Surveys  \vol 35 \yr 1980 \pages 1--80 \endref


\ref\key St \by A. Starkov \book Dynamical systems in \hs s 
\bookinfo FAZIS, Moscow,
1999; English translation in:  Translations of Mathematical Monographs, vol. 190
\publ Amer. Math. Soc.
\publaddr  Providence, R.I. \yr 2000    \endref

\ref\key U1 \by  M. Urba\'nski \paper  Diophantine
approximation of self-conformal measures\paperinfo Preprint \yr  2003 \endref Ê

\ref\key U2 \bysame \paper  Diophantine
approximation for conformal measures of one-dimensional iterated
function systems\paperinfo Preprint \yr  2003 \endref Ê



\ref\key W \by B. Weiss  \paper  Almost no points on a
Cantor set are very well approximable \jour Proc. R. Soc. Lond. \vol 
457\yr 2001 \pages  949--952 \endref 

\ref \key Z \by F. \v Zeludevi\v c \paper Simultane diophantische
Approximationen abh\"angiger Gr\"ossen in mehreren Metriken \jour Acta
Arith.
\vol  46
\yr 1986 \issue 3 \pages 285--296 \lang German \endref

\endRefs
\enddocument

\end